\documentclass[a4paper,12pt]{amsart}
\usepackage{amsfonts, amssymb, euscript,graphics}

\usepackage{amsmath}
\usepackage{amscd}
\usepackage{eepic}
\usepackage{comment}
\usepackage{psfrag}
\usepackage{wrapfig}
\newcommand{\R}{{\mathbb R}}

\newcommand{\Z}{{\mathbb Z}}

\renewcommand{\le}{\leqslant}
\renewcommand{\ge}{\geqslant}

\renewcommand{\em}[1]{{\par\noindent\bf #1}}

\newtheorem{statement}{Statement}[section]
\newfont{\cmbb}{cmex10 at 7pt}

\newcommand{\ign}[1]{{}}
\def \.{\mskip 1mu}
\def \?{\mskip -1mu}

\renewcommand*\proof[1]{{\noindent\it Proof#1.\/\enspace\ignorespaces}}
\newcommand*\pr{\proof{}}

\def\sq{\raisebox{-.1ex}{\hbox{$\square$}}}

\def\qed{{\unskip\nobreak\hfill\hskip.5em\nobreak\hfil\sq
\parfillskip=0pt\medbreak}}

\makeatletter

\renewcommand{\abovecaptionskip}{.8ex}

\newcommand*{\capdelimiter}{.}
\renewcommand{\@makecaption}[2]{%
\vspace{\abovecaptionskip}%
\sbox{\@tempboxa}{#1\capdelimiter#2}
\ifdim \wd\@tempboxa >\hsize
   #1\capdelimiter#2 \par
\else
 \global\@minipagefalse \hbox to \hsize {\hfil #1\capdelimiter#2 \hfil}%
\fi
\vspace{\belowcaptionskip}}

\newcommand{\subs}[1]
   {\refstepcounter{subsection}
   \medskip\noindent
   {\it\arabic{section}.\arabic{subsection}.\hspace{.25em}\ignorespaces#1.}
          \ignorespaces}

   \@addtoreset{equation}{section}

\title{On Hamiltonian and contact isotopy liftings}

\author{Petr E. Pushkar\\\today}
\address[Petr E. Pushkar]{National Research University Higher School of Economics (HSE)}
\email{petya.pushkar@gmail.com}
\thanks{The work  were partially supported by the grant RFBR 15-01-05990.}

\begin{document}


\begin{abstract}
We construct counterexamples to lifting properties of Hamiltonian and contact isotopies.
\end{abstract}
\maketitle
\date{AAAA}

\section*{Introduction}  The main result of this  paper is a construction of counterexamples to lifting of isotopies in symplectic and contact settings. 
In the Introduction we give definitions of symplectic reduction and isotopy lifting, the first section  of the paper is devoted to the contact part. 
The last section contains counterexamples to the lifting property.

\subs{Symplectic reduction along a hypersurface} We recall the notion of a symplectic reduction along a hypersurface \cite{AG}. 
A smooth hypersurface in a symplectic manifold carries a field of characteristic directions which are kernels of the restriction of the symplectic structure to the hypersurface. 
We will suppose that the space of characteristics (i.e. integral curves of that field) is a manifold (it always holds locally) and, in particular, that the natural projection of the hypersurface to the space of characteristics is a smooth fibration. 
The manifold of characteristics has a natural symplectic structure descendig from the initial symplectic structure \cite{AG}. 
The resulting symplectic manifold is called a symplectic reduction along the hypersurface.

Let $M$ be a symplectic reduction along a hypersurface $\widetilde{Z}$ in a symplectic manifold $\widetilde{M}$. Denote by $\widetilde{\Pi}$ the natural projection $\widetilde{Z} \to M$. 
For a subset $\widetilde{L}$ in $\widetilde{M}$ its symplectic reduction along $\widetilde{Z}$ is, by definition, $\widetilde{\Pi}(\widetilde{L}\cap\widetilde{Z})$. 
It is well known that for a Lagrangian submanifold in $\widetilde{M}$, which is transversal to $\widetilde{Z}$, its symplectic reduction along $\widetilde{Z}$ is an immersed Lagrangian submanifold in $M$.  

It is said that an isotopy $I(t)$ $(t\in [0,1])$ of $L$ ($L$ is a symplectic reduction of $\widetilde{L}$) lifts to $\widetilde{M}$ if there exists an isotopy $\widetilde{I}(t)$ $(t\in [0,1])$ of  $\widetilde{M}$ such that $I(t)(L)$ is a symplectic reduction of $\widetilde{I}(t)(\widetilde{L})$.

\subs{Lifting of an isotopy} The following Hamiltonian isotopy lifting property was claimed for the case of general (i.e. not only along a hypersurface) symplectic reduction in \cite{EG} (Lemma 2.5.1):

\begin{statement} 
\label{sympiso}
If the subset $\widetilde{L} \subset \widetilde{M}$ is closed and the projection of $\widetilde{L}\cap \widetilde{Z}$ to M
is proper, then every compact Hamiltonian isotopy $I(t)$ in $M$ lifts to a compact Hamiltonian isotopy $\widetilde{I}(t)$ in $\widetilde{M}$ which maps some neighbourhood $\widetilde{Z_0} \subset \widetilde{Z}$ of $\widetilde{L}\cap \widetilde{Z}\subset \widetilde{Z}$ into $\widetilde{Z}$ for all $t \in [0,1]$.
\end{statement}

If the hypersurface $\widetilde{Z}$ is compact and cooriented then it is easy to prove this statement. 
Any smooth compactly supported extension $\widetilde{H}_t$ of a peimage of a Hamiltonian $H_t$ of the isotopy on $M$ gives us a Hamiltonian for the required isotopy on $\widetilde{M}$. 
Surprisingly, in contrast with an analogues statement in the smooth category, that proof could not be adapted to the general non-compact situation and the statement is false in general. 
Formal mistake in the proof of \ref{sympiso} in \cite{EG} is an assumption that the Hamiltonian vector field of $\widetilde{H}_t$ is globally integrable (i.e. gives rise to a phase flow defined for all values of time) on $\widetilde{Z}$. 
We show below that in general one needs to have  a ``sufficient room'' in $\widetilde{M}$ to construct a lifting of a compactly supported isotopy of a reduction of a Lagrangian submanifold.

\subs{Simplest case}  The simplest and the most ``visible'' case illustrating our construction and techniques lies in the contact topology setting, and we are going to describe it without giving definitions  of contact reduction and contact isotopy lifting. 
Consider the space $J^1(\R)=J^1(\R,\R)$ of 1-jets of functions on the real line. 
$J^1(\R)$ is a contact manifold with the contact structure given by the form $du-pdq$ (where $q$ is the coordinate on the source, $u$ on the target, $p$ is the derivative).
The $1$-jet extension $j^1(f)$ of a function $f$ on $\R$ is a Legendrian manifold in $J^1(\R)$, here $j^1(f)=\{(p,q,u)\in J^1(\R)| p=\frac{df}{dq} (q), u=f(q)\}$. 
Consider a Legendrian isotopy $\Lambda_{t, t\in [0,1]}$ of the $1$-jets extension of zero function ($\Lambda_0=j^1(0)$). 
Let also $\Lambda_{t, t\in [0,1]}$ be supported in an open subset $A$ of $J^1(\R)$ given by inequalities $|q|<1$, $|p|<~1$.
We claim that $\Lambda_1$ must contain a point with zero $q$-coordinate, such that the absolute value of its $u$ coordinate is at most~$1$. 
Indeed, the front (the projection on the $(u,q)$-plane) of $\Lambda_1$ contains a graph of a continuous function $h$ (it is a well-known corollary from Chekanov theorem \cite{Cha,Ch}). 
After a small perturbation one can assume that $h$ is piece-wise smooth  whose  derivative is (absolutely) bounded by 1 in its smooth points, since the isotopy is supported in $A$ where $|p|<1$. Since $h$ is supported in $[-1,1]$ then, by mean value theorem, $C^0$-norm of $h$ is bounded by one and thus $|h(0)|<1$.
Denote by $u(x)$ the value of $u$-coordinate of the point $x$. We can conclude that not for every smooth function $f$ on $[0,1]$ (such that $f(0)=0$) there exists an isotopy $\Lambda_{t,t\in[0,1]}$ of $\Lambda_0=j^1(0)$ supported in $A$ such that $u(\Lambda_t\cap \{q=0\})=\{f(t)\}$ for every $t\in [0,1]$.  
It means, as we will see below, that not every legendrian isotopy (of a point in a 1-dimensional contact manifold in the example described) admits a lifting. 
Similar phenomenon exists in the multydimensional case as well as in the symplectic case and is described below.

\section{Contact reduction along a hypersurface and isotopy lifting} We recall the notion of contact reduction along a hypersurface. The following construction is a contact analog of  symplectic reduction along a hypersurface above. 
Consider a contact manifold $(\widetilde{N},\widetilde{\xi})$ and a smooth hypersurface $\widetilde{Z}$ in $\widetilde{N}$. 
We will suppose that at each point of $\widetilde{Z}$ the tangent plane to  $\widetilde{Z}$ is transversal to $\xi$. 
In that case $\xi$ carries a natural (characteristic) field of directions (\cite{AG}). 
Similarly to the symplectic case we suppose that its integral curves form a manifold $N$, so that the natural projection $\Pi \colon \widetilde{Z} \to N$ is a smooth fibration. 
Then $N$ carries a canonical contact structure $\xi$ which is uniquely defined by the condition, that the preimage of $\xi$ under the action of $d\Pi$ concides with 
$\widetilde{\xi} \cap T\widetilde{Z}$.

Definitions of  contact reduction of a subset along a hypersurface and lifting of an isotopy coincide with the symplectic case above. 
The following contact isotopy lifting property was formulated (modulo few misprints) in \cite{EG} (Lemma 2.6.1):

\begin{statement} 
\label{contiso}
If the subset $\widetilde{L}$ is closed in $\widetilde{N}$ and the projection of $\widetilde{L}\cap \widetilde{Z}$ to $N$
is proper, then for every compact contact isotopy $I(t), t \in [0,1]$ in $N$ there exists a compact contact isotopy $\widetilde{I}(t)$ in $\widetilde{N}$ 
such that \\
(a) $\Pi (\widetilde{I}(t)\widetilde{L} \cap \widetilde{Z}) = I(t)L$,\\
(b) $\widetilde{I}(t)$ maps some neighbourhood of $\widetilde{L} \cap \widetilde{Z} \subset \widetilde{Z}$ into $\widetilde{Z}$ for all $t \in [0,1]$.
\end{statement} 

\section{Counterexamples}
We will show that Statements  \ref{contiso} and \ref{contiso} are not true in general. 

\subs{Contact case} Let us start from the contact case and describe $\widetilde{N}, \widetilde{Z}$ and $\widetilde{L}$. Let $M$ be a closed connected manifold.
Consider the contact manifold $J^1(M\times \R)=T^*M\times T^*\R\times\R$, denote by $(q,p)$ the canonical coordinates on the factor $T^*\R$.
Let $\widetilde{N}$ be an open submanifold of $J^1(M\times \R)$ given by inequalities  $|q|<1,|p|<1$. 
$\widetilde{N}$ is a contact manifold with the induced contact structure. 
We define $\widetilde{Z}$ to be a hypersurface in $\widetilde{N}$ given by an equation $q=0$.
Consider the 1-jet extension $j^1\widetilde{f}\subset J^1(M\times \R)$ of zero function $\widetilde{f}$ on $M\times\R$, 
let $\widetilde{L}$ be $j^1\widetilde{f}\cap\widetilde{N}$.
For such a hypersurface $\widetilde{Z}$ the contact reduction $N$ is naturally contactomorphic to the space $J^1(M)$ with its standard contact structure, $L$ is the
1-jet extension $j^1f$ of zero function $f$ on $M$. 
The hypersurface $\widetilde{Z}$ and Legendrian manifold $j^1f$ satisfy conditions of Statement \ref{contiso}.
Note that for any smooth function $g$ on $M$ $j^1g$ could be connected to $j^1f$ by a compactly supported contact isotopy. 

\begin{statement}
\label{ineq}
Consider a compactly supported contact isotopy $\varphi_{t, t\in [0,1]}$ of the manifold $J^1(M)$.  
Suppose that $\varphi_1(j^1f)=j^1g$ for a function $g\colon M \to \R$. 
Suppose there exists a compactly supported lifting of that isotopy to a contact isotopy on $\widetilde{N}$ satisfying conditions of Statement \ref{contiso}. 
Then $$\max\limits_{x \in M}g(x)\le 1, \max\limits_{x \in M}g(x) - \min\limits_{x \in M}g(x) \le 2$$. 
\end{statement} 

As a corollary we get that it is impossible to lift any contact isotopy joining $j^1f$ with $j^1g$ where $\max\limits_{x \in M}g(x)>1$. 
That result is non-trivial even if $M$ is a point, in fact it is exactly the case that has been described in the Introduction.

\pr Consider a compactly supported contact isotopy of $\widetilde{N}$ and extend it trivially to the compactly supported contact isotopy 
$\widetilde{\varphi}_{t, t\in [0,1]}$ (supported in $\widetilde{N}$) of  $J^1(M\times\R)$. 
By Chekanov theorem (\cite{Ch}) there exists a $K$ such that $\widetilde{\varphi}_1(j^1\widetilde{f})$ is given by a quadratic at infinity (with respect to $\R^K$) generating family $F\colon M\times\R\times\R^K \to \R$. 
Consider now the function $F$ as a family $G_{q, q\in\R}$ of functions depending on $q\in\R$, $\R$ is the factor in $M\times\R\times\R^K$, $G_q(x,w)=F(x,q,w)$ for $x \in M, w\in \R^K$.
Recall that a Cerf diagramm of a family of functions is a graph of all critical values.
Consider the Cerf diagramm $\Gamma$ of ${G_q}$. For $|q|\ge 1$  all critical values of  $G_q(x,w)$ are equal to zero.
Cerf diagramm $\Gamma$ contains a graph of continuous function $h\colon [-1,1] \to \R$, such that 
$h(0)=\max\limits_{x \in M}g(x)$.
Existence of such a function follows from Viterbo's theory of selected values of quadratic at infinity functions. We briefly sketch it here.

For a quadratic at infinity function $G\colon M\times\R^K \to \R$ and non-zero homological class $\alpha\in H_*(M,\Z_2)$ one can correspond a critical value $c(\alpha, G)$ of function $G$ in a following way. 
Denote by $G^a$ the sublevel set $\{x \in M\times \R | G(x)\le a\}$.
For a sufficiently big number $C>0$ the pair $(G^C,G^{-C})$ is naturally homotopy equivalent to the Thom space of a trivial bundle over $M$.
Denote by $T$ the Thom isomorphism $H_*(M;\Z_2)\to H_*(G^C,G^{-C};\Z_2)$.  
Then $c(\alpha,G)= \inf \{a|T\alpha \in i_* H_*(G^a,G^{-C};\Z_2) \}$, where $i$ is a natural inclusion map $(G^a,G^{-C})\to (G^C,G^{-C})$.
It turns out that $c(\alpha,G)$ is a critical value of $G$ and it depends on $G$ continuously.
Consider a function $h\colon \R \to \R$, $h(q)=c([M],G_q)$, $[M]$ is $\Z_2$-fundamental class of $M$.
Consider the number $h(0)$. Note that $G_0$ is a quadratic at infinity generating family for $j^1g$. 
We claim that if $G$ is a quadratic at infinity generating family for 1-jet extension $j^1g$ of a function $g\colon M\to \R$ then $c([M], G)=\max\limits_{x \in M}g(x)$.
The fastest way to prove it is to use Theret uniqueness theorem \cite{Th}, claiming that, in that particular case, that $G$ is equivalent, after stabilization, to a quadratic stabilization of $g$. 
For stabilizations (i.e. functions of type $g(x)+Q(w)$, $Q$ is nondegenerate quadratic form) that statement is obvious.

The graph of the function $h$ is, by definition, a subset in $\Gamma$. 
Note that $h(q)=0$ for $|q|\ge 1$, since all the critical values of $G_q$ are zero.  
For generic small perturbation  $\widetilde{\varphi}^*_{t, t\in [0,1]}$ of $\widetilde{\varphi}_{t, t\in [0,1]}$ supported in $\widetilde{N}$ the corresponding perturbed family $G^*_{q, q\in\R}$ 
is sufficiently generic and the corresponding perturbed function $h^*$ is piece-wise smooth. 
For any $t\in [0,1]$ Legendrian manifold $\widetilde{\varphi}_t(j^1\widetilde{f})$ is transverse to $\widetilde{Z}$. 
Hence, the reduction along $\widetilde{Z}$ of $\widetilde{\varphi}_1(j^1\widetilde{f})$ is an 1-jet extension of a function $g^*$ which is $C_0$-close to $g$.
So the number $h^*(0)$ is close to $\max\limits_{x \in M}g(x)$.
For a generic point $q_0\in ]-1,1[$ (except a countable discreet set in $]-1,1[$) in a neighborhood $U(q_0)$ of $q_0$ the critical point corresponding to $c([M], G^*_q)$ smoothly depends on $q$: $h^*(q)=F^*(x(q),q,w(q))$ for smooth functions $x(q),w(q)$, such that for any $q\in U(q_0)$ $(x(q),w(q))$ is a critical point of $G^*_q$.
The point $w(q)$ is a critical point of $F^*(x(q),q,.)\colon \R^N \to \R$ and generate a point $l(q) \in \widetilde{\varphi}^*_1(j^1\widetilde{f})$.
The absolute value of the $p$-coordinate of the point $l(q)$ is at most $1$, since $l(q)\in \widetilde{N}$.
By the chain rule the derivative $\frac{dh^*}{dq}(q)$ equals to that  $p$-coordinate. Thus $h^*(0) < 1$, since $h^*$ is a piece-wise smooth continuous function whose derivative is at least 
$-1$ and $h^*(1)=0$. 

Since $\max\limits_{x \in M}g(x)$ is close to $h^*(0)$ we get the first inequality of the statement. Similarly, using the class  $[pt]\in H_0(M;\Z_2)$, we get 
$\min\limits_{x \in M}g(x)\ge~-1$. Hence,   $\max\limits_{x \in M}g(x) - \min\limits_{x \in M}g(x) \le 2$.
\qed

\subs{Symplectic case} A counterexample in the symplectic case is similar to the contact case above. Let $M$ be a closed connected manifold of  dimension at least 1.
Consider the symplectic manifold $T^*M\times T^*\R$, denote by $(q,p)$ the canonical coordinates on the factor $T^*\R$.
Let $\widetilde{N}$ be an open submanifold of $T^*M\times T^*\R$ given by inequalities  $|q|<1,|p|<~1$. 
$\widetilde{N}$ is a symplectic manifold with the induced contact structure. 
We define $\widetilde{Z}$ to be a hypersurface in $\widetilde{N}$ given by the equation $q=0$.
Consider the graph $\Gamma(d\widetilde{f})\subset T^*M\times T^*\R$ of differential of the zero function $\widetilde{f}$ on $M\times\R$, 
let $\widetilde{L}$ be $\Gamma(d\widetilde{f})\cap\widetilde{N}$.
For such a hypersurface $\widetilde{Z}$ the symplectic reduction $N$ is naturally symplectomorphic to the space $T^*M$ with its standard symplectic structure, 
$L$ is the graph $\Gamma(df)\subset T^*M\times T^*\R$ of zero function on $M$. 
The hypersurface $\widetilde{Z}$ and Lagrangian manifold $\Gamma(df)$ satisfy conditions of Statement \ref{sympiso}.

\begin{statement}
Consider a compactly supported Hamiltonian isotopy $\psi_{t, t\in [0,1]}$ of the manifold $T^*M$.  
Suppose that $\psi_1(\Gamma(df))=\Gamma(dg)$ for a function $g\colon M \to \R$. 
Suppose there exists a compactly supported lifting of that isotopy to a Hamiltonian isotopy on $\widetilde{N}$ satisfying conditions of Statement \ref{contiso}. 
Then $\max\limits_{x \in M}g(x) - \min\limits_{x \in M}g(x) \le 2$. 
\end{statement} 

\pr Consider a lifting isotopy  $\widetilde{\psi}_t$ on $T^*M\times T^*\R$. 
We get a family $\{\widetilde{\psi}_t(\Gamma(d\widetilde{f}))\}$ of Lagrangian maifolds. 
Each of them is an exact Lagrangian manifold, thus we can cover that isotopy by compactly supported isotopy of 
Legendrian manifolds $\widetilde{\Lambda_t}\subset J^1(M\times \R)$ projecting to  $\{\widetilde{\psi}_t(\Gamma(d\widetilde{f}))\}$ and 
coinciding with $j^1\widetilde{f}$ outside $\widetilde{N}\times\R$.
The reduction of $\widetilde{\Lambda}_1$ along $\widetilde{Z}\times \R$ is $j^1(g')$, where $g'$ differs from $g$ by a constant.
To the family $\Lambda_t$ we can apply considerations of Statement  \ref{ineq} . Concluding  inequality on $g'$ obviously finish the proof.  
\qed

Any two graphs of differentials of functions could be joined by compactly supported Hamiltonian isotopy.  Hence it is impossible to lift a Hamiltonian isotopy joining the 
graph of the differential of zero function with the graph of differential of a function $g$ such that $\max\limits_{x \in M}g(x) - \min\limits_{x \in M}g(x) > 2$. 
Such a function $g$ obviously exists if the dimension of $M$ is not zero.

\end{document}